\documentclass[letterpaper,oneside,10pt]{article}
\usepackage[USenglish]{babel} %francais, polish, spanish, ...
\usepackage[ansinew]{inputenc}
\usepackage{verbatim}
\usepackage{amsmath, amsthm, amssymb, amsfonts}
\usepackage[USenglish]{babel} %francais, polish, spanish, ...
\usepackage{verbatim}
\usepackage[pdftex]{graphicx}
\usepackage{color}
\usepackage{float}
\usepackage{latexsym}
\usepackage{amsmath}
\usepackage{amsfonts}
\usepackage{enumerate}
\usepackage{url}
\usepackage{algorithm2e}
\begin{document}
\newcommand{\sg}{\mathbf{s}_g}
\newcommand{\bg}{\mathbf{b}_g}
\newcommand{\sa}{\mathbf{s}_a}
\newcommand{\ba}{\mathbf{b}_a}
\newcommand{\sgh}{\hat{\mathbf{s}}_g}
\newcommand{\bgh}{\hat{\mathbf{b}}_g}
\newcommand{\sah}{\hat{\mathbf{s}}_a}
\newcommand{\bah}{\hat{\mathbf{b}}_a}
\newcommand{\mw}{\left\|\omega\right\|}
\newcommand{\dd}{\partial}
%%\pagestyle{empty} %No headings for the first pages.

%% Title Page %%%%%%%%%%%%%%%%%%%%%%%%%%%%%%%%%%%%%%%%%%%%%%%
\pagestyle{plain}

\title{\textbf{Online Performance Optimization of a DC Motor Driving a Variable Pitch Propeller}}

\author{Raphael Cohen \footnote{French graduate school ENSMA, ``Ecole Nationale Sup\'{e}rieur de M\'{e}canique et d'A\'{e}rotechnique", Poitiers, France, {\small cohen.raphael79@yahoo.com.}},
David Miculescu \footnote{School of Aerospace Engineering, Massachusetts Institute of Technology, Cambridge, MA {\small dmicul@mit.edu.}},
Kevin Reilley \footnote{School of Aerospace Engineering, Georgia Institute of Technology, Atlanta, GA 30332, {\small kreilley3@gatech.edu.}},
Mehrdad Pakmehr \footnote{postdoctoral fellow at the School of Aerospace Engineering, Georgia Institute of Technology, Atlanta, GA 30332, {\small mehrdad.pakmehr@gatech.edu.}},
Eric Feron \footnote{professor at the School of Aerospace Engineering, Georgia Institute of Technology, Atlanta, GA 30332, {\small feron@gatech.edu.}}      }
%\end{center}

\date{\null}

\maketitle
\pagestyle{plain} %Now display headings: headings / fancy / .
%%%%%%%%%%%%%%%%%%%%%%%%%%%%%%%%%%%%%%%%%%%%%%%%%%%%%%%%%%%%%%%%%%%%%%

\section*{Abstract}

A practical online optimization scheme is developed for performance optimization of an electrical aircraft propulsion system. The goal is to minimize the power extraction of the propulsion system for any given thrust value. The online optimizer computes the optimum pitch angle of a variable pitch propeller by minimizing the power of the system for a command thrust value. This algorithm is tested on a DC motor driving a variable pitch propeller; the experimental hardware setup of the DC motor along with its variable pitch propeller is also described. Experimental results show the efficiency and practicality of the proposed online optimization scheme. Outstanding issues are sketched.

%%%%%%%%%%%%%%%%%%%%%%%%%%%%%%%%%%%%%%%%%%%%%%%%%%%%%%%%%%%%%%%%%%%%%%
%%%%%%%%%%%%%%%%%%%%%%%%%%%%%%%%%%%%%%%%%%%%%%%%%%%%%%%%%%%%%%%%%%%%%%%%%%%%%%%%
\section{Introduction}
%%%%%%%%%%%%%%%%%%%%%%%%%%%%%%%%%%%%%%%%%%%%%%%%%%%%%%%%%%%%%%%%%%%%%%%%%%%%%%%%%%%%%%%%%%%%%%%%%%%%%%%%%%%%%%%%%%%%%%%%%
Variable pitch propellers are used in various systems including ships, airplanes, helicopters, and other dynamical systems. In aerospace applications, it is well-known that variable pitch propellers are used to adapt the system to different thrust levels and air speeds so that the propeller blades don't stall; on the other hand, if the pitch angle is not optimized, the propulsive efficiency will be reduced. To keep the system performance at its optimum level, there is a need for continuous online optimization algorithms for the propulsion systems driving variable pitch propellers. Some of the recent research works regarding the control and optimization of systems with variable pitch propellers are described in \cite{ModContProp-dullens-2009, ContPropOpt-balsamo-2011, cpp-Martelli-2013, PeakSeekProp-cazenave-2011, cazenave-msthesis-2012, windTurbine-biegel-2011, windturbine-lianyou-2011, pitchcontrol-muljadi-2001}. The examples of variable pitch propeller optimization for fast ferries and ships are presented in \cite{ModContProp-dullens-2009, ContPropOpt-balsamo-2011, cpp-Martelli-2013}, and for aerospace applications are presented in \cite{PeakSeekProp-cazenave-2011, cazenave-msthesis-2012}. Variable pitch systems for windmill optimization problems can be found in \cite{windTurbine-biegel-2011, windturbine-lianyou-2011, pitchcontrol-muljadi-2001}.

Some of the noteworthy research on the optimization for aerospace systems are mentioned in the following. Traditional helicopters rely on single RPM rotor and modulate rotor thrust by varying pitch. A-160 Hummingbird helicopter \cite{a160-link} uses variable RPM to optimize autonomy of vehicle. A method for controlling the thrust of an aircraft engine using a single lever power controller is described in \cite{patent-vos-2002}; the engine control commands comprise propeller RPM and engine inlet manifold air pressure commands.  The present paper builds upon these prior works by proposing an approach for online optimization of propulsion system efficiency. Moreover, our work demonstrates the feasibility of the approach by describing an experimental implementation of the proposed approach.

This work seeks to implement a code which optimizes online the power in a DC motor with a variable pitch propeller. The DC motor represents more complex propulsion systems with aerospace applications. However, it is valuable in its own right, given that many small-UAV propulsion systems are powered by DC motors, and an increasing number of large air vehicles rely on electric engines as well \cite{helios-link, solarimpulse-link}. The goal is to minimize the power extraction of an aero propulsion system for any given thrust value. This is important since propulsion systems may operate in different operation conditions. For aerospace systems, these variations happen by changes in the altitude and also the weather conditions and speed. Another type of variation happens in the propulsion hardware due to aging, damage, and maintenance. So, there is a need for online optimization algorithms which help the system to always operate with the best possible performance.

A DC motor with a variable pitch propeller with two blades has been constructed for experiments at zero air speed. The optimization of the system for non-zero air speeds is also studied via simulation. Then, a code that optimizes the power for a given level of thrust is developed and implemented. Experiments show the efficiency and practicality of the proposed online optimization scheme for simple propulsion systems like DC motors with variable pitch props. This method could be extended and implemented to more complex aero engines such as turboprops and turboshafts with variable pitch actuation.  The method could also be used in turbofan applications if secondary actuation schemes, such as variable fan or compressor inlet vanes are available.  The rest of this paper is organized as follows. In section 2, we go through the modeling and present the equations of motion for the system. In section 3, we discuss the experimental setup. In section 4, we justify the need for online optimization by performing computer simulations. In section 5, we present the online optimization scheme for a DC motor driving a variable pitch propeller. In section 6, we conclude the paper.

%%%%%%%%%%%%%%%%%%%%%%%%%%%%%%%%%%%%%%%%%%%%%%%%%%%%%%%%%%%%%%%%%%%%%%%%%%%%%%%%%%%%%%%%%%%%%%%%%%%%%%%%%%%%%%%%%%%%%%%%%
\section{Problem Formulation}
%%%%%%%%%%%%%%%%%%%%%%%%%%%%%%%%%%%%%%%%%%%%%%%%%%%%%%%%%%%%%%%%%%%%%%%%%%%%%%%%%%%%%%%%%%%%%%%%%%%%%%%%%%%%%%%%%%%%%%%%%
We intend to develop a real-time power optimization algorithm for a DC motor driving a variable pitch propeller. The following model describes a system comprising an electric DC motor and a variable pitch propeller. This is the elementary representation of propulsion systems with a prop/fan, and is useful for pedagogical purposes illustrating a simple optimization algorithm for these types of systems. As aircraft and rotorcraft systems increasingly rely on electric propulsion systems, the model is becoming increasingly relevant to real-world applications as well.

\begin{figure}[!ht]
\centering
\includegraphics[width=0.45\textwidth]{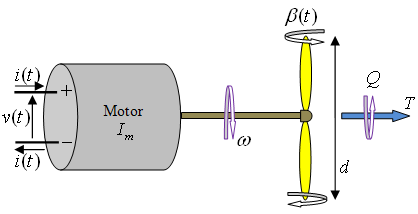}
\caption{Schematic of a DC motor driving a variable pitch propeller }
\label{fig:MotorPropScheme}
\end{figure}

\begin{figure}[!ht]
\centering
\includegraphics[width=0.3\textwidth]{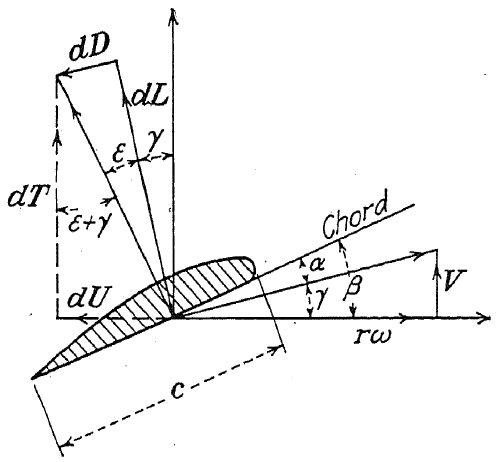}
\caption{Resolution of forces acting on a propeller section \cite{mises-theory}}
\label{fig:PropForces}
\end{figure}

The equations of motion for a DC motor driving a variable pitch propeller system are
\begin{equation}\label{eqn_1}
\begin{array}{l}
\displaystyle \frac{d}{dt}\omega(t) = \frac{1}{I_m}\left[k_b i(t)-B_1 \omega(t)-Q(\omega(t),\beta(t)) \right],\\[5pt]
\displaystyle \frac{d}{dt}i(t) = \frac{1}{L}\left[ -R i(t)-k_b \omega(t) + v(t) \right],\\[5pt]
\displaystyle y(t)=T(\omega(t),\beta(t)),
\end{array}
\end{equation}
where  $v(t)$: voltage across the armatures (V), $i(t)$: current through the motor (A),  $\omega(t)$  motor shaft angular speed (rad/s),  $I_m$: motor inertia $(kg.m^2/s^2)$, $k_b$: induced emf constant (V.s/rad), $B_1$: viscous friction coefficient (N.m.s/rad), $R$: resistance (Ohm), $L$: self-inductance (H), $Q$: propeller torque (N.m), $T$: propeller thrust (N), and $\beta(t)$: propeller pitch angle (rad). The schematic illustration of a DC motor and a variable pitch propeller is shown in Figure \ref{fig:MotorPropScheme}.
For computing the forces acting on a propeller blade, we use \textit{blade element theory}. Figure \ref{fig:PropForces} shows a cross section of the blade. The forces which the air exerts on the blade element of the spanwise extension $dr$ is the resultant of two component forces, lift and drag, perpendicular and parallel, respectively, to the velocity with which the element moves through the air.
The velocity appears in Figure \ref{fig:PropForces} as a vector with the horizontal component $r\omega=2\pi r n$ and the vertical component is $V$. The square of the resulting velocity is
\begin{equation}\label{eqn_2}
\displaystyle V^2+(r\omega)^2 = V^2+(2\pi r n)^2 = n^2 d^2 \left[ J^2 + (2\pi r /d)^2 \right],
\end{equation}
where $J=(V/nd)$ denotes the advance ratio, $n$ denotes the round per second (RPS), and $d$ denotes the propeller diameter. Let $G:=J^2 + (2\pi r /d)^2$. The angle between the velocity vector and the horizontal axis is $\gamma = \arctan (V/r\omega)$. If $\beta$ denotes the angle between the plane of rotation and blade section chord, the angle of attack is $\alpha = \beta - \gamma$.

In the following it is assumed that, for all the blade sections, $C_L = C_{L_{\alpha}}(\alpha - \alpha_0)$ and $C_D = C_{D_0}+ k C_{L_{\alpha}}^2(\alpha - \alpha_0)^2$ are functions of the angle $\alpha$. The lift and drag forces for the blade element between the radii $r$ and $r+dr$ can be obtained as:
\begin{equation}\label{eqn_3}
\begin{array}{l}
\displaystyle dL=0.5 C_L \rho V^2 dS \\[5pt]
\displaystyle ~~~~ =\frac{1}{2} C_{L_{\alpha}}(\alpha - \alpha_0) \rho n^2 d ^2 G c(r) dr,\\[5pt]
\displaystyle dD=0.5 C_D \rho V^2 dS \\[5pt]
\displaystyle ~~~~ =\frac{1}{2} (C_{D_0}+ k C_{L_{\alpha}}^2(\alpha - \alpha_0)^2) \rho n^2 d ^2 G c(r) dr,
\end{array}
\end{equation}
where $c(r)$ denotes the chord length of the blade section at the radii $r$. It will be useful to decompose the resultant of the lift and drag forces into the components $dT$ and $dU$ parallel and perpendicular to the propeller axis. These components are given by
\begin{equation}\label{eqn_4}
\begin{array}{l}
\displaystyle dT=dL \cos \gamma - dD \sin \gamma,\\[5pt]
\displaystyle dU=dL \sin \gamma + dD \cos \gamma.
\end{array}
\end{equation}
The overall forces on the propeller are
\begin{equation}\label{eqn_thrust}
\begin{array}{l}
\displaystyle T=0.5 m \rho n^2 d^2 \int^{d/2}_{0} G (C_L \cos \gamma - C_D \sin \gamma) c(r) dr,
\end{array}
\end{equation}
\begin{equation}\label{eqn_U}
\begin{array}{l}
\displaystyle U=0.5 m \rho n^2 d^2 \int^{d/2}_{0} G (C_L \sin \gamma + C_D \cos \gamma) c(r) dr,
\end{array}
\end{equation}
where $m$ is the number of blades of the propeller. $dU$ has the moment $r dU$ with respect to the propeller axis. It contributes to the propeller torque $Q$ and to power $P$, which can be written as
\begin{equation}\label{eqn_torque}
\begin{array}{l}
\displaystyle Q=m \int^{d/2}_{0} r dU \\[5pt]
\displaystyle ~~~ = 0.5 m \rho n^2 d^2 \int^{d/2}_{0} G (C_L \sin \gamma + C_D \cos \gamma) r c(r) dr,
\end{array}
\end{equation}
\begin{equation}\label{eqn_power}
\begin{array}{l}
\displaystyle P= \omega Q = 2\pi nm \int^{d/2}_{0} r dU \\[8pt]
\displaystyle ~~~ = \pi nm \rho n^2 d^2 \int^{d/2}_{0} G (C_L \sin \gamma + C_D \cos \gamma) r c(r) dr.
\end{array}
\end{equation}
In this work, we assume $\epsilon=\alpha$. As it can be gathered from (\ref{eqn_thrust}) and (\ref{eqn_power}), thrust and power are functions of parameters $n$, $V$, and $\beta$. Its is not easy to compute thrust and power using these equations; hence, there is a need to come up with a numerical approach to compute them.

%%%%%%%%%%%%%%%%%%%%%%%%%%%%%%%%%%%%%%%%%%%%%%%%%%%%%%%%%%%%%%%%%%%%%%%%%%%%%%%%%%%%%%%%%%%%%%%%%%%%%%%%%%%%%%%%%%%%%%%%%
\section{Experimental Hardware Setup}
\begin{figure}[!ht]
\centering
\includegraphics[width=0.5\textwidth]{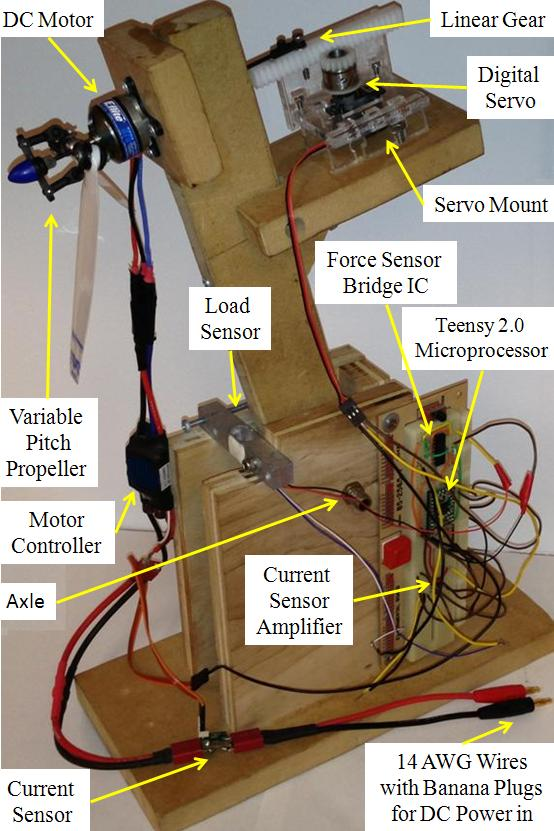}
\caption{Experimental hardware setup of the DC motor with its variable pitch propeller } \label{fig:Experimental_hardware_setup1}
\end{figure}

%%%%%%%%%%%%%%%%%%%%%%%%%%%%%%%%%%%%%%%%%%%%%%%%%%%%%%%%%%%%%%%%%%%%%%%%%%%%%%%%%%%%%%%%%%%%%%%%%%%%%%%%%%%%%%%%%%%%%%%%%
In this section, a description of the experimental setup is given. The components of the setup include a DC motor, a variable pitch propeller, a motor controller, a control microprocessor, a pitch angle actuator, a current-voltage sensor, and a load sensor (strain gauge). Figure \ref{fig:Experimental_hardware_setup1} shows the experimental setup.

The DC motor chosen for this experiment is an \textit{E-Flite Park 370 BL Outrunner 1200Kv} \cite{Reference3} motor. This motor has been designed for electric models equipped with variable pitch propeller systems, and it is powerful enough to generate a non-negligible thrust (about 10 Oz (300 g)). This DC motor, which drives a variable pitch propeller, can be commanded by an Electronic Speed Control (ESC) unit; the ESC that is used in this work is an ALIGN RCE-BL35X ESC~\cite{Reference4}. The power supply plugs into the DC motor through the ESC, and the ESC is itself controlled by a Teensy 2.0 microprocessor \cite{Reference8}. In other words, the voltage sent from the microprocessor to the ESC will affect the power going to the DC motor. Thanks to the presence of the \textit{Arduino} system, the DC motor can be commanded by a computer program. \textit{Arduino} is the software which can be used to load programs on the microprocessor. All the programs in Arduino can be coded using C++.

\begin{figure}[!ht]
\centering
\includegraphics[width=0.6\textwidth]{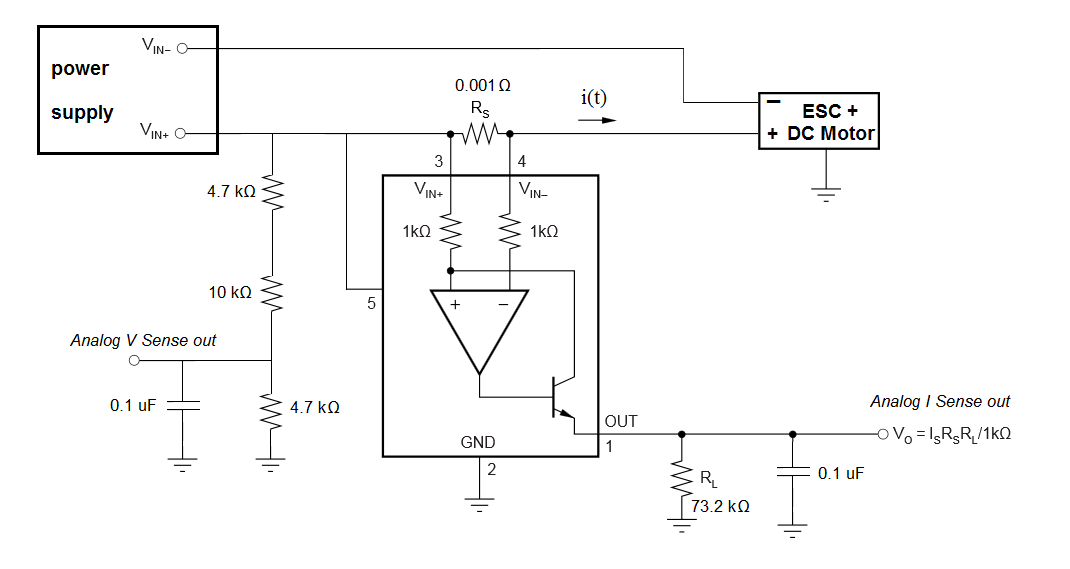}
\caption{Electric diagram of the IV Sensor ~\cite{Reference6}} \label{fig:Electric Diagram IVSensor}
\end{figure}

In order to measure the propeller thrust, a strain gauge is installed on the system. This strain gauge is supposed to have an electrical resistance proportional to its deformation. Due to the linear behavior of the load sensor, the voltage across the load sensor is proportional to the thrust. As the resting voltage of the sensor is low relative to the threshold voltages necessary for the microcontroller analog-to-digital converter, a single $INA125U$ \cite{Reference9} operational amplifier was inserted into the circuit between the load sensor and the microcontroller.  To measure the motor power, a current-voltage sensor (ATTOPILOT 45A~\cite{Reference6}) is used. Figure \ref{fig:Electric Diagram IVSensor} shows the electrical connections between the power supply, sensor, ESC, and the DC motor.
With this sensor, both current and voltage can be calculated by measuring an output voltage.
\begin{equation}\label{eqn:output_voltage}
\displaystyle v_{an, v}=\frac{4.7}{4.7+10}v(t)=0.32v(t),
\end{equation}
\begin{equation}\label{eqn:output_current}
\displaystyle v_{an, i}=\frac{R_{s} R_{L}}{1 k}i(t)= 0.073i(t).
\end{equation}
where, $v_{an, v}$ and $v_{an, i}$ are the output voltages from the IV Sensor that are proportional to the voltage and current going into the DC motor. $R_{s}$ and $R_{L}$ are resistors; the electric diagram in Figure \ref{fig:Electric Diagram IVSensor} shows these resistors.  The voltage in the DC motor varies between $7~V$ and $12~V$, so the output voltage is between $2.24~V$ and $3.84~V$, whereas the current varies between $0$ and $1~A$, and the output current voltage is between 0 and $73~mV$, which is very low. So, for practical purposes, the output voltage from the current sensor is amplified. The amplifier used is the $LM386N$ amplifier \cite{Reference10}.

\subsection{Calibration of the Sensors} %%%%%%%%%%%%%%%%%%%%%%%%%%%%%%
To calibrate the voltage and current sensors, a linear regression is used. Both sensors are well adapted and the output values are affine functions of the values that needed to be measured. These calibrations are valid for any voltage and for currents higher than 300 mA threshold. The relations of these calibrations are as follows: $v=0.0202 v_{an}-0.0237 ~V$, $i=4.1532 i_{an}-1826.67 ~mA$, where $v_{an}$ and $i_{an}$ are the analog voltages theoretically proportional to the voltage and current, $v$ and $i$, going into the DC motor. For the load sensor, it is known that the output voltage is proportional to the thrust. The voltage indicated for zero thrust and the voltage related to a known load (that of a hanging object whose weight is known) are available; as a result the calibration relation is $T= 9.81((4.11v_{out}-49.36)/1000))~(N)$.

\begin{figure}[!ht]
\centering
\includegraphics[width=0.6\textwidth]{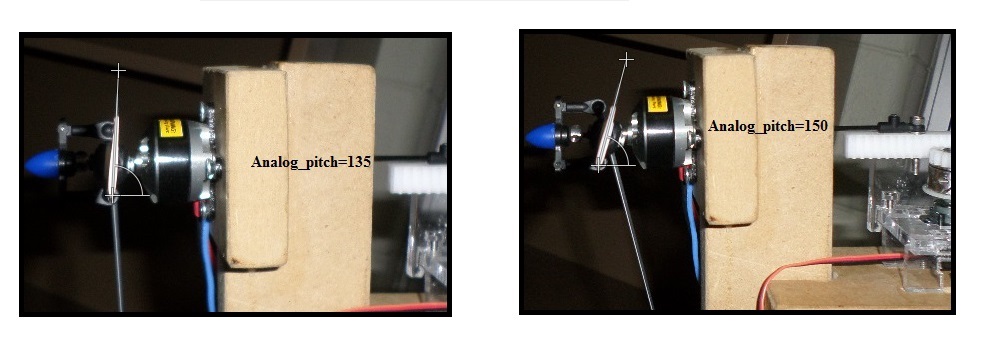}
\caption{Pitch angle calibration using GIMP} \label{fig:pitch_angle_calib}
\end{figure}

It is also important to calibrate the pitch angle in order to have an idea of the blades' angle. The pitch angle is controlled by a servo, and the only parameter we know about the pitch angle is the input voltage of this servo. So, for different values of the servo's input voltage, some pictures of the propeller are taken, and then, the angle of the blades is measured using Figure \ref{fig:pitch_angle_calib} an the \textit{GIMP} software \cite{gimp-link}. The calibration relation is $\beta=0.59(v_{servo}-135)~(deg)$.

%%%%%%%%%%%%%%%%%%%%%%%%%%%%%%%%%%%%%%%%%%%%%%%%%%%%%%%%%%%%%%%%%%%%%%%%%%%%%%%%%%%%%%%%%%%%%%%%%%%%%%%%%%%%%%%%%%%%%%%%%
\section{Numerical Simulations}
%%%%%%%%%%%%%%%%%%%%%%%%%%%%%%%%%%%%%%%%%%%%%%%%%%%%%%%%%%%%%%%%%%%%%%%%%%%%%%%%%%%%%%%%%%%%%%%%%%%%%%%%%%%%%%%%%%%%%%%%%
In order to justify the need for online optimization, some numerical simulation codes are now described. The codes simulate the DC motor dynamics driving a variable pitch propeller. The results of the codes are the numerical computation of the power and thrust for variable pitch angle values. Numerical simulations show the convexity of the functions, an important property, which makes optimization considerably easier.

\subsection{Matlab Code: $P=f(\beta)$ Plot}%%%%%%%%%%%%%%%%%%%%%%%%%%%%%%%%

To investigate the possibility of DC motor power optimization for non-zero air speed, motor power versus pitch angle curves are plotted for given values of thrust and air speed. Air speed, $V$, and thrust, $T$, are known and we would like to calculate the necessary DC motor power for different values of pitch angle.
\begin{figure}[!ht]
\centering
\includegraphics[width=0.55\textwidth]{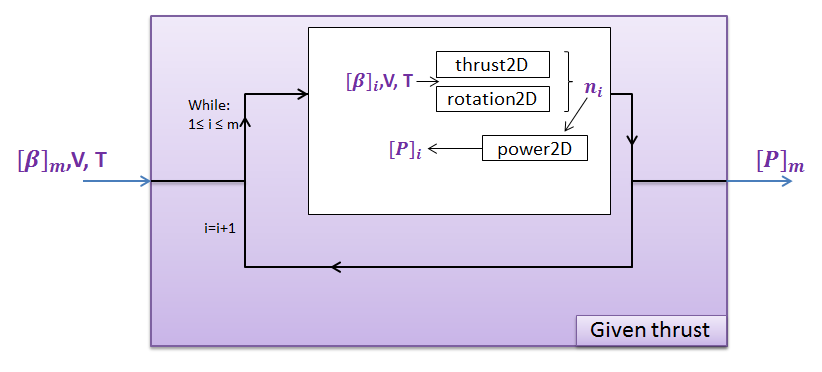}
\caption{Structure of the Matlab code to plot $P=f(\beta)$} \label{fig:structure MATLAB 2D}
\end{figure}
As it is known from equation (\ref{eqn_thrust}), the propeller thrust depends on $V$, $\beta$, and $n$. So, if $\beta$, $V$, and $T$ are known, we would be able to find $n$. For the case where $V=0$, and $\beta$ is known, equation (\ref{eqn_thrust}) can be solved easily to obtain $n$ (equation of the type $n^3=k$). But, when $V \ne 0$, the parameter $n$ appears in the computation of the integral, yielding an implicit equation for $n$. To overcome this problem, a code with three main functions is developed that calculates the points $(\beta,P(\beta))$. Figure \ref{fig:structure MATLAB 2D} is a schematic visualization of this Matlab code, which plots the power versus pitch angle for a given thrust. The three functions of this code are as follows:
\begin{enumerate}
\item \textit{thrust2D}: this function calculates the thrust by integration given $n$, V and $\beta$ . $n$ is now known and the integration is possible.
\item \textit{rotation2D}: this function calculates the thrust given the same $V$ and $\beta$ but for any rotational speed $n$, and generates a polynomial approximation of the thrust, $T=f(n)$, using parameters $\beta$ and $V$. Using this approximation, we can find the rotational speed, $n$, corresponding to the original given level of thrust, $T_c$.
\item \textit{power2D}: this function calculates the power of the DC motor given $T_c$, $\beta$, and $V$ using (\ref{eqn_power}). All the parameters are now known thanks to the \emph{rotation2D} and \emph{thrust2D} functions.
\end{enumerate}

\begin{figure}[!ht]
\centering
\begin{minipage}[b]{0.48\textwidth}
\centering
\includegraphics[width=0.99\textwidth]{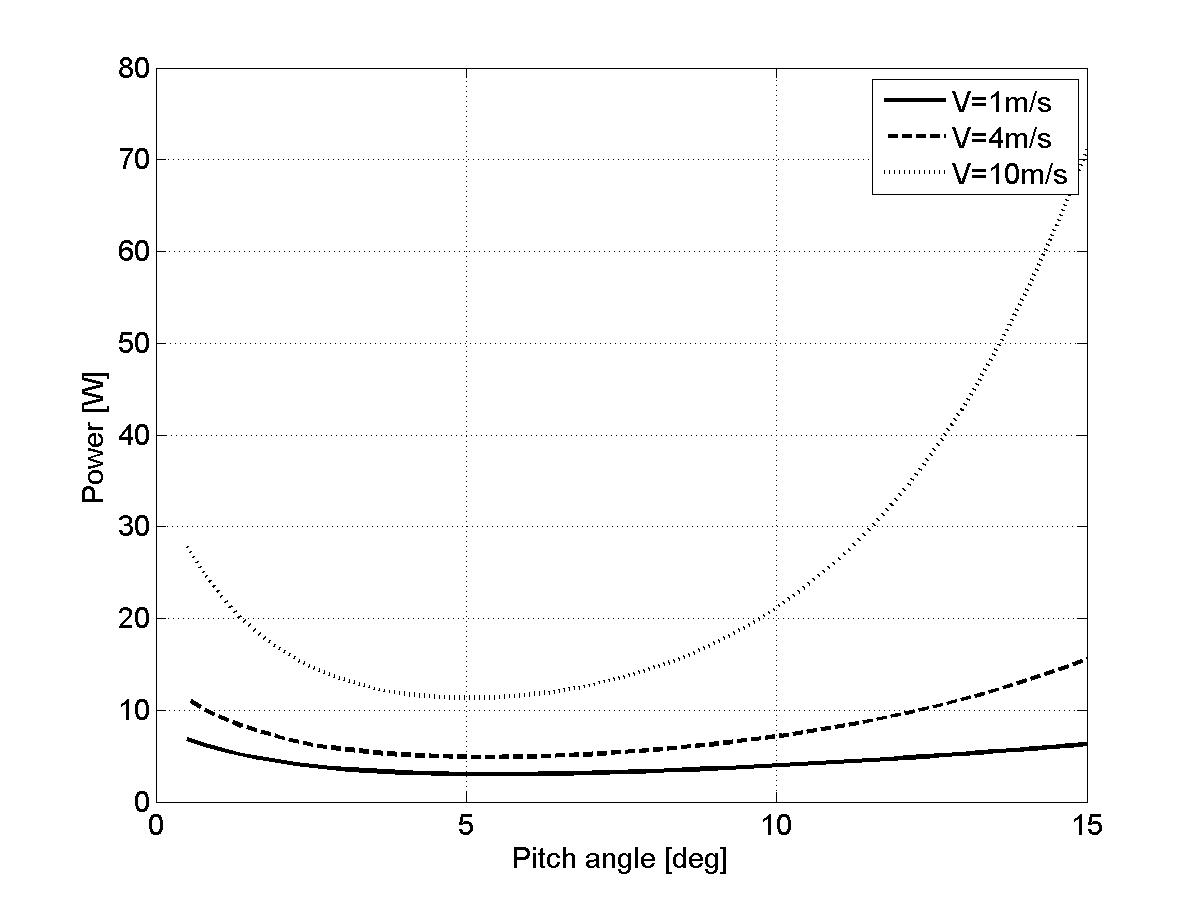}
\caption{Power vs pitch angle for different air speeds, and $T_c=0.3~(N)$  } \label{fig:2D_different_AirSpeeds}
\end{minipage}
\hfill
\begin{minipage}[b]{0.48\textwidth}
\centering
\includegraphics[width=0.99\textwidth]{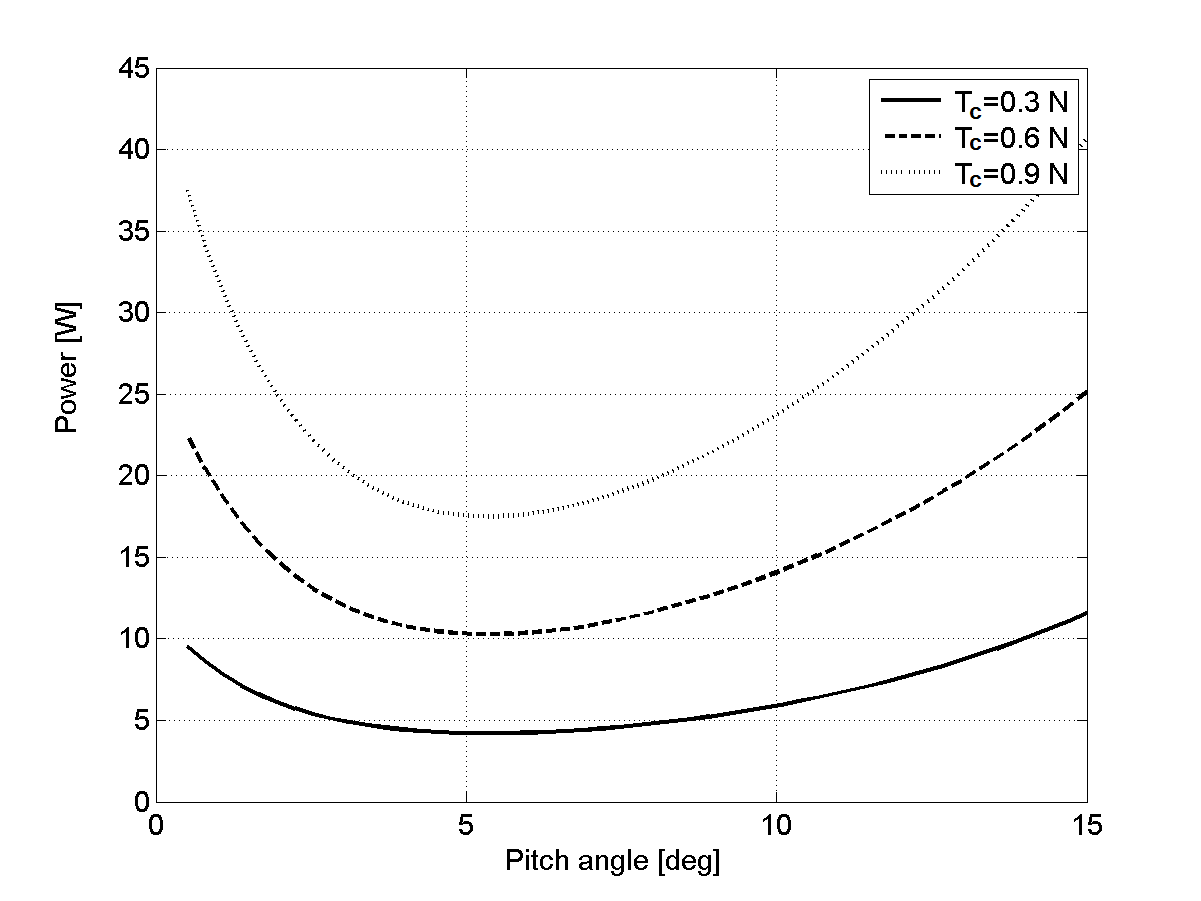}
\caption{Power vs pitch angle for different thrust commands, and $V=3~(m/s)$ } \label{fig:2D_different_thrusts}
\end{minipage}
\end{figure}

To understand what is the impact of the air speed and the thrust command variations on the optimization process, simulations are performed. Figure \ref{fig:2D_different_AirSpeeds} shows the curves of power vs pitch angle for three different air speeds, $V=1,~4,~10~(m/s)$. For each of these three simulations, it can be observed that a single optimal angle exists; the optimal angle slowly increases with the increase in the air speed. This observation reinforces the idea of the development of an online optimizer on the pitch angle during the aircraft flight. Figure \ref{fig:2D_different_thrusts} shows curves of power vs pitch angle for different thrust commands $T_c=0.3,~0.6,~0.9~(N)$. As it can be observed, for various thrust commands, the optimal power does not seem to be very sensitive to blade pitch angle, suggesting that a single blade pitch angle is satisfactory for more applications. However, as will be seen later, experimental measures show very different behaviors.

\subsection{Matlab Code: $T=f(P,\beta)$ Plot}%%%%%%%%%%%%%%%%%%%%%%%%%%%%%%%%%%%%%%%%%%%%

\begin{figure}[!ht]
\centering
\includegraphics[width=0.55\textwidth]{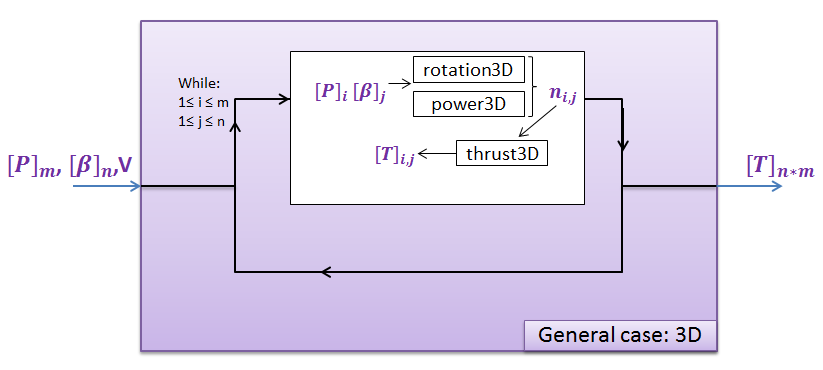}
\caption{Structure of the Matlab code calculating the thrust for given pitch angle and power} \label{fig:structure MATLAB 3D}
\end{figure}
\begin{figure}[!ht]
\centering
\includegraphics[width=0.64\textwidth]{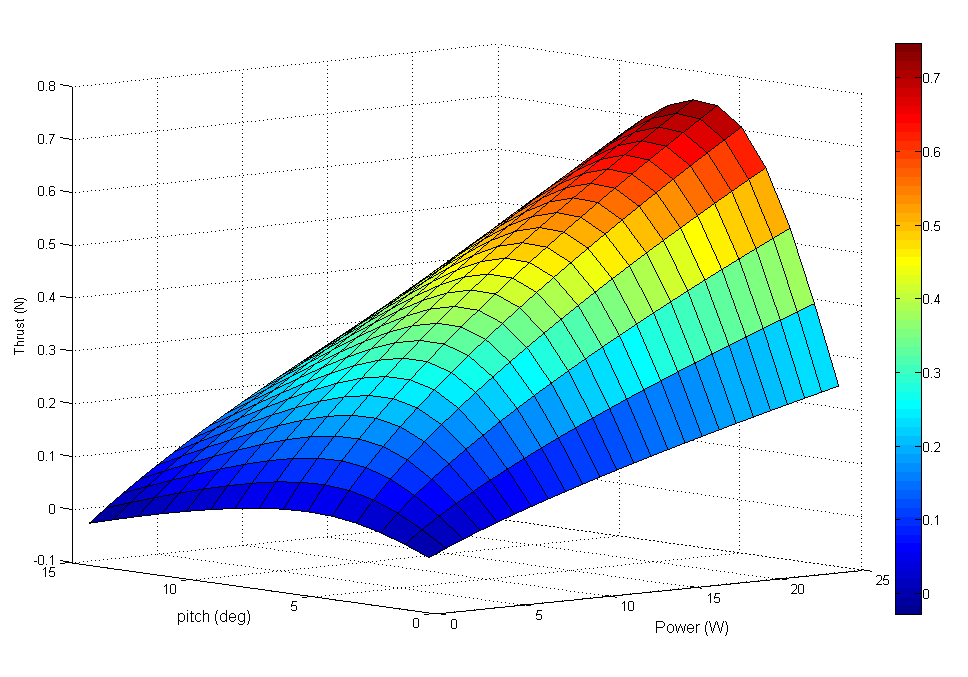}
\caption{3D visualization of thrust versus power and pitch angle for $V=3~(m/s)$} \label{fig:MATLAB 3D black and white}
\end{figure}

To calculate the propeller thrust for given pitch angle and motor power, another code is developed. Figure \ref{fig:structure MATLAB 3D} shows the structure of the Matlab code calculating the thrust versus power and pitch angle. The code developed here is similar to the previous approach with the only exception that equation (\ref{eqn_power}) is also used in calculation process. For given values of power and pitch angle, equation (\ref{eqn_power}) can be used to obtain rotational speed ($n$). The three functions of the code are described below:
\begin{enumerate}
\item \textit{power3D}: this function calculates motor power for given values of $\beta$, $V$ and $n$.
\item \textit{rotation3D}: this function calls function \textit{power3D}, for given values of $\beta$, $V$, and for several values of rotational speed, $n$. Then, it uses polynomial approximation approach to compute power ($P=f(n)$) for the given values of $\beta$ and $V$. Finally, thanks to the polynomial approximation, this function returns the rotational speed $n$ corresponding to a given $P$ and $\beta$.
\item \textit{thrust3D}: this function computes the thrust for given values of $\beta$, $V$, and $n$.
\end{enumerate}
Figure \ref{fig:MATLAB 3D black and white} shows a plot of the thrust versus motor power and propeller pitch angle. As it is apparent from this plot, for any value of thrust, power can be minimized by finding an optimum pitch angle (i.e., thrust is a convex function of motor power and blade pitch).

%%%%%%%%%%%%%%%%%%%%%%%%%%%%%%%%%%%%%%%%%%%%%%%%%%%%%%%%%%%%%%%%%%%%%%%%%%%%%%%%%%%%%%%%%%%%%%%%%%%%%%%%%%%%%%%%%%%%%%%%%
\section{Online Performance Optimization Scheme}
%%%%%%%%%%%%%%%%%%%%%%%%%%%%%%%%%%%%%%%%%%%%%%%%%%%%%%%%%%%%%%%%%%%%%%%%%%%%%%%%%%%%%%%%%%%%%%%%%%%%%%%%%%%%%%%%%%%%%%%%%

\subsection{Control with Fixed Pitch Input}
\begin{figure}[!ht]
\centering
\includegraphics[width=0.8\textwidth]{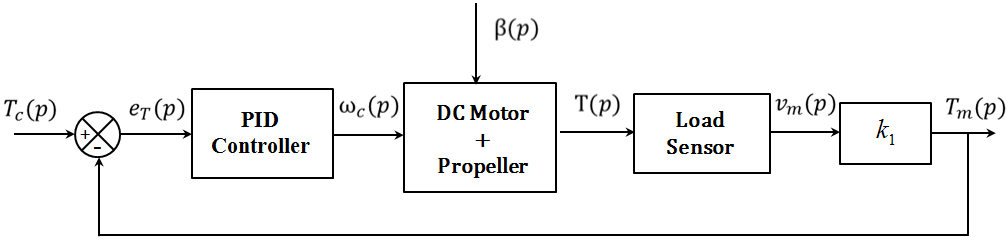}
\caption{Block Diagram of the system with a PID controller} \label{fig:PID_block diagram}
\end{figure}

\begin{figure}[!ht]
\centering
\begin{minipage}[b]{0.48\textwidth}
\centering
\includegraphics[width=0.99\textwidth]{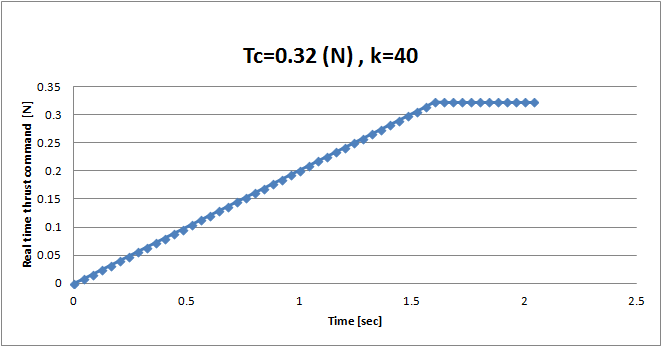}
\caption{Commanded thrust ramp from 0 to $T_c=0.32~(N)$ in 1.6 sec} \label{fig:thrust_40_points}
\end{minipage}
\hfill
\begin{minipage}[b]{0.48\textwidth}
\centering
\includegraphics[width=0.99\textwidth]{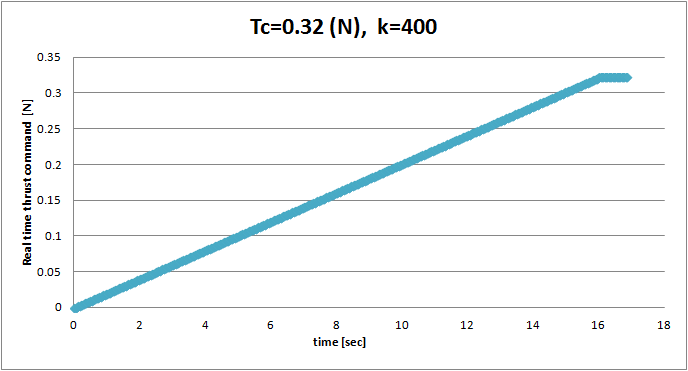}
\caption{Commanded thrust ramp from 0 to $T_c=0.32~(N)$ in 16 sec} \label{fig:thrust_400_points}
\end{minipage}
\end{figure}

Here, the behaviour of the system with speed controller is investigated for fixed values of the propeller pitch command. Figure \ref{fig:PID_block diagram} shows a block diagram of the system with a PID controller for rotational speed control. In this diagram, $\omega_{c}(p)$ is the command rotational speed; $v_{m}(p)$ is the measured voltage of the load sensor; $\beta(p)$ is the fixed pitch angle of the blades; $T_c(p)$ is the thrust command; $T_m(p)$ is the measured thrust; $e_T(p)$ is the error in the thrust; and $k_1$ is an affine computation block that converts load sensor voltage ($v_{m}(p)$), into a thrust value ($T_m(p)$). The PID controller computes a rotational speed command for the DC motor for a reference thrust value, while the pitch angle is fixed.
The performance of the control system is evaluated by feeding saturated ramps for $T_c(p)$. Indeed, it was observed that the response of the system to step functions was inappropriate, maybe due to non-optimal PID controller gain tuning and strong nonlinear effects on the closed-loop system response. PID tuning was not pursued any further than necessary for our investigation, since power optimization is the centerpiece of this paper. However, further optimization of the engine PID controller is still performed at the time of the writing of this paper. Figures \ref{fig:thrust_40_points} and \ref{fig:thrust_400_points} show two examples of saturated ramp commands: Both ramps saturate at the same thrust level ($0.32 N$). However, the first ramp (shown in Fig. \ref{fig:thrust_40_points}) is steeper than the second (shown in Fig. \ref{fig:thrust_400_points}).

\begin{figure}[!ht]
\centering
\begin{minipage}[b]{0.48\textwidth}
\centering
\includegraphics[width=0.99\textwidth]{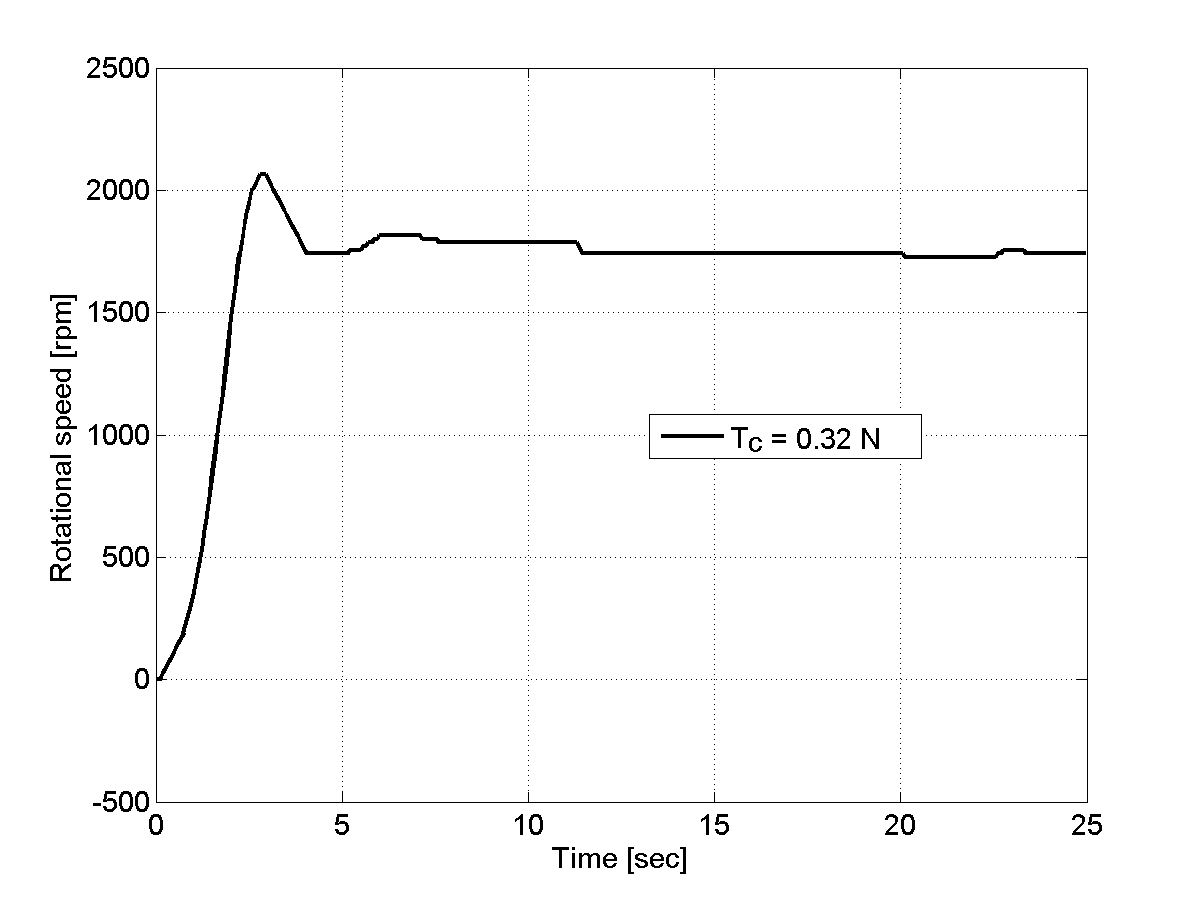}
\caption{Response (in RPM) to first ramp (Fig \ref{fig:thrust_40_points})} \label{fig:pid_response_40_points}
\end{minipage}
\hfill
\begin{minipage}[b]{0.48\textwidth}
\centering
\includegraphics[width=0.99\textwidth]{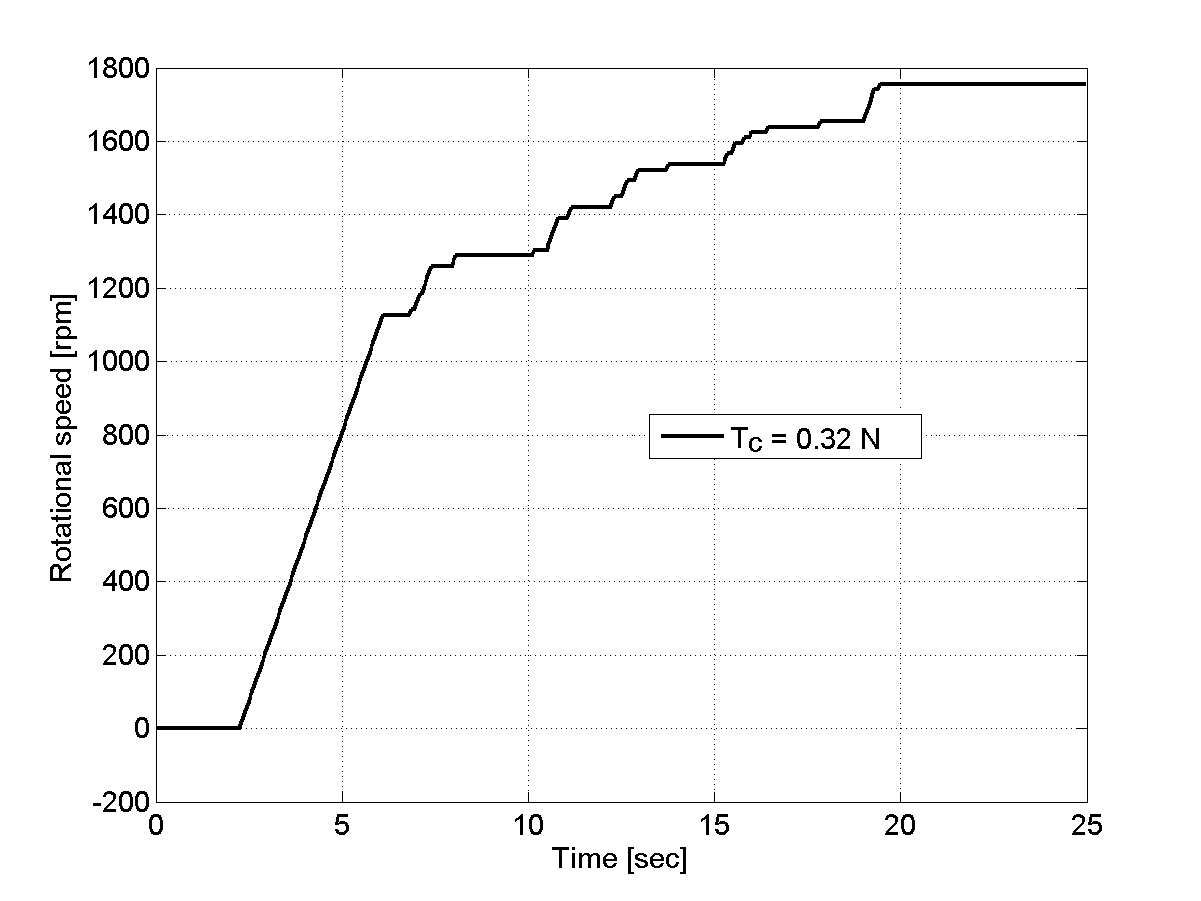}
\caption{Response (in RPM) to second ramp (Fig \ref{fig:thrust_400_points}} \label{fig:pid_response_400_points}
\end{minipage}
\end{figure}

Figures \ref{fig:pid_response_40_points} and \ref{fig:pid_response_400_points} show the time history of the rotational speed for a given level of thrust for for the two ramps discussed above. The response to the first ramp features a $12.5\%$, while the response to the second ramp is too slow. In what follows, saturated ramps where the ramp lasts 1.6 sec are used.

\subsection{Basic Optimization: Fixed Pitch Steps}%%%%%%%%%%%%%%%%%%%%%%%%%%%%%%
\begin{figure}[!ht]
\centering
\includegraphics[width=0.5\textwidth]{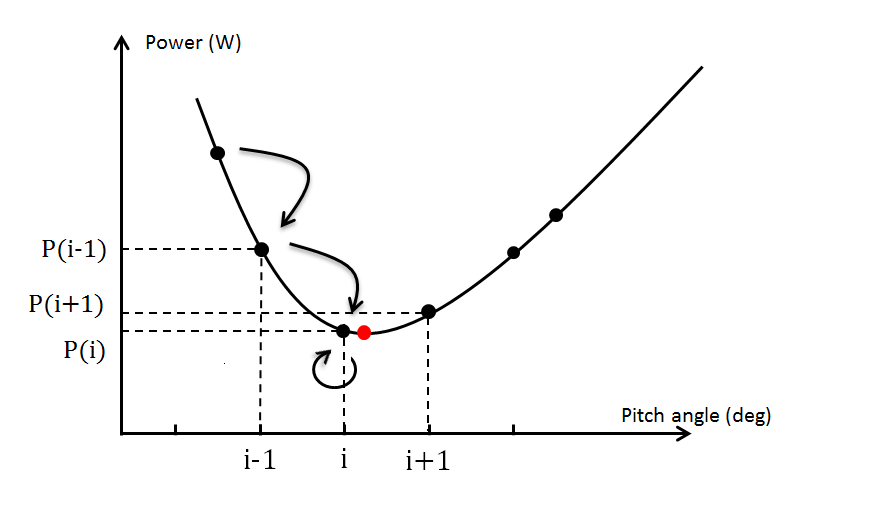}
\caption{Functional diagram for optimization process} \label{fig:schema_optimization_1}
\end{figure}

\begin{figure}[!ht]
\centering
\includegraphics[width=0.65\textwidth]{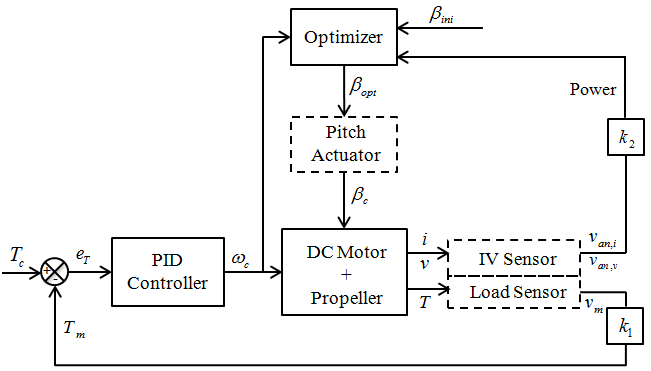}
\caption{Schematic of the closed-loop system with its online optimizer} \label{fig:Combined_PID_Optimizer}
\end{figure}

The first experiment is performed using a simple, fixed-step optimization scheme. It is already known that for a given level of thrust, power versus pitch angle curve is a convex curve, (see Figures \ref{fig:2D_different_AirSpeeds} and \ref{fig:2D_different_thrusts}). Based on this fact, an optimization function is developed that calculates the required motor power to achieve a given thrust level. Depending on whether required power decreases or increases, the next next step is chosen to be positive or negative. Figure \ref{fig:schema_optimization_1} shows the functional diagram of the optimization process and Figure \ref{fig:Combined_PID_Optimizer} show a diagram of the closed-loop system with the online optimizer.

The fixed step length optimization process is presented as a pseudocode in Algorithm \ref{Alg:FixedPitch}.
\begin{algorithm}[!ht]
 \SetAlgoLined
 {\small
 initialization\;
 $\beta=0.59~deg$;
 \newline
 $T_c=0.52~N$;
 \newline
 direction=1;\newline
 pitch\_step=0.59;\newline
 set\_propeller ($\beta,T_c$); \newline
 measured prev\_power; \newline
 diff\_power=1; \newline
 \While{3 minutes}{
  \eIf{power saturation}{
   direction=1;}{\eIf{diff\_power$>0$}{direction=1;}{direction=-1;}}
   $\beta=\beta$ + direction $\times$ pitch\_step; \newline
   set\_propeller ($\beta,T_c$); \newline
   measure next\_power; \newline
   diff\_power=prev\_power-next\_power; \newline
   prev\_power=next\_power;     }  }
 \caption{Pseudocode of the basic optimization approach with fixed pitch steps}\label{Alg:FixedPitch}
\end{algorithm}
In Algorithm \ref{Alg:FixedPitch}, function \textit{set\_propeller} calls the program with the PID controller that makes the system to reach the thrust command for a given pitch angle $\beta$.

\begin{figure}[!ht]
\centering
\begin{minipage}[b]{0.48\textwidth}
\centering
\includegraphics[width=0.99\textwidth]{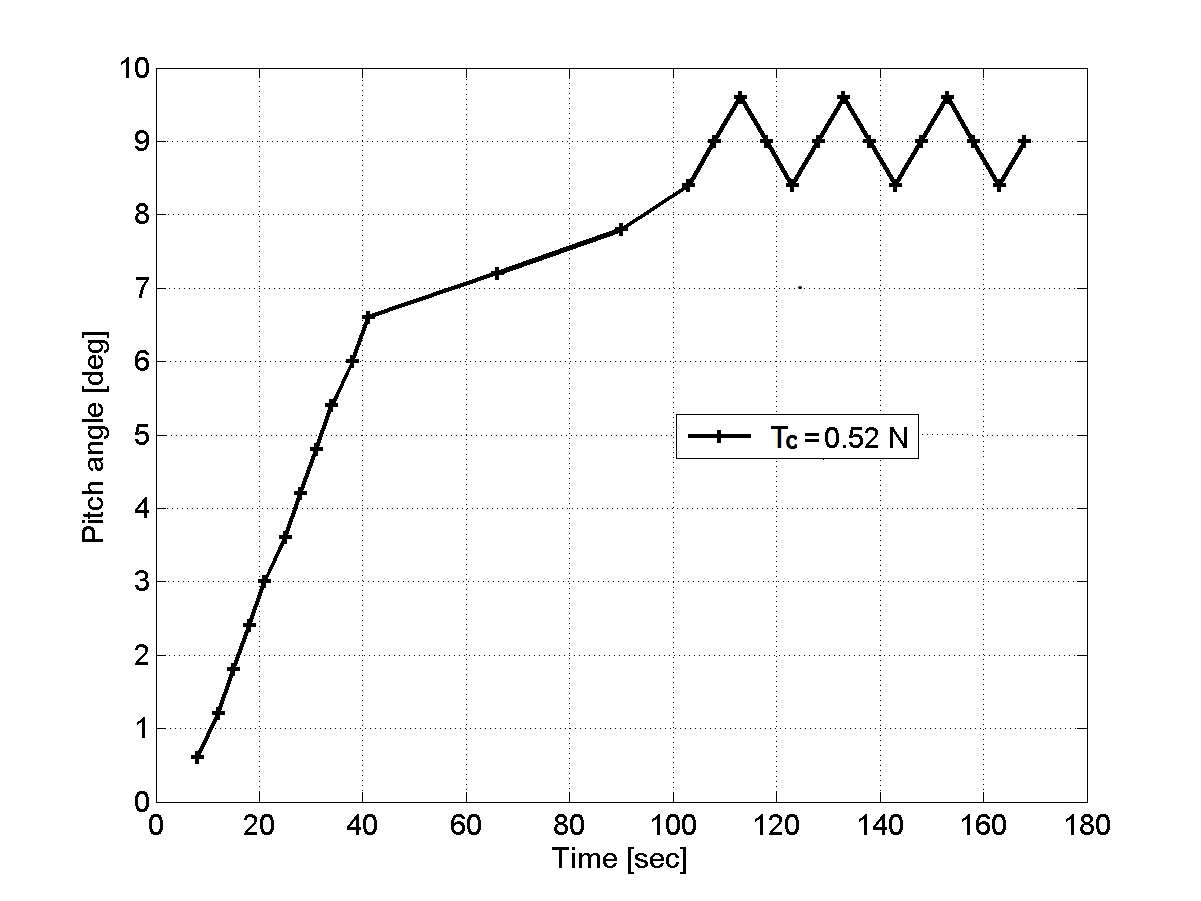}
\caption{Pitch angle of the propeller during the online optimization process} \label{fig:pitch_optimization_3}
\end{minipage}
\hfill
\begin{minipage}[b]{0.48\textwidth}
\centering
\includegraphics[width=0.99\textwidth]{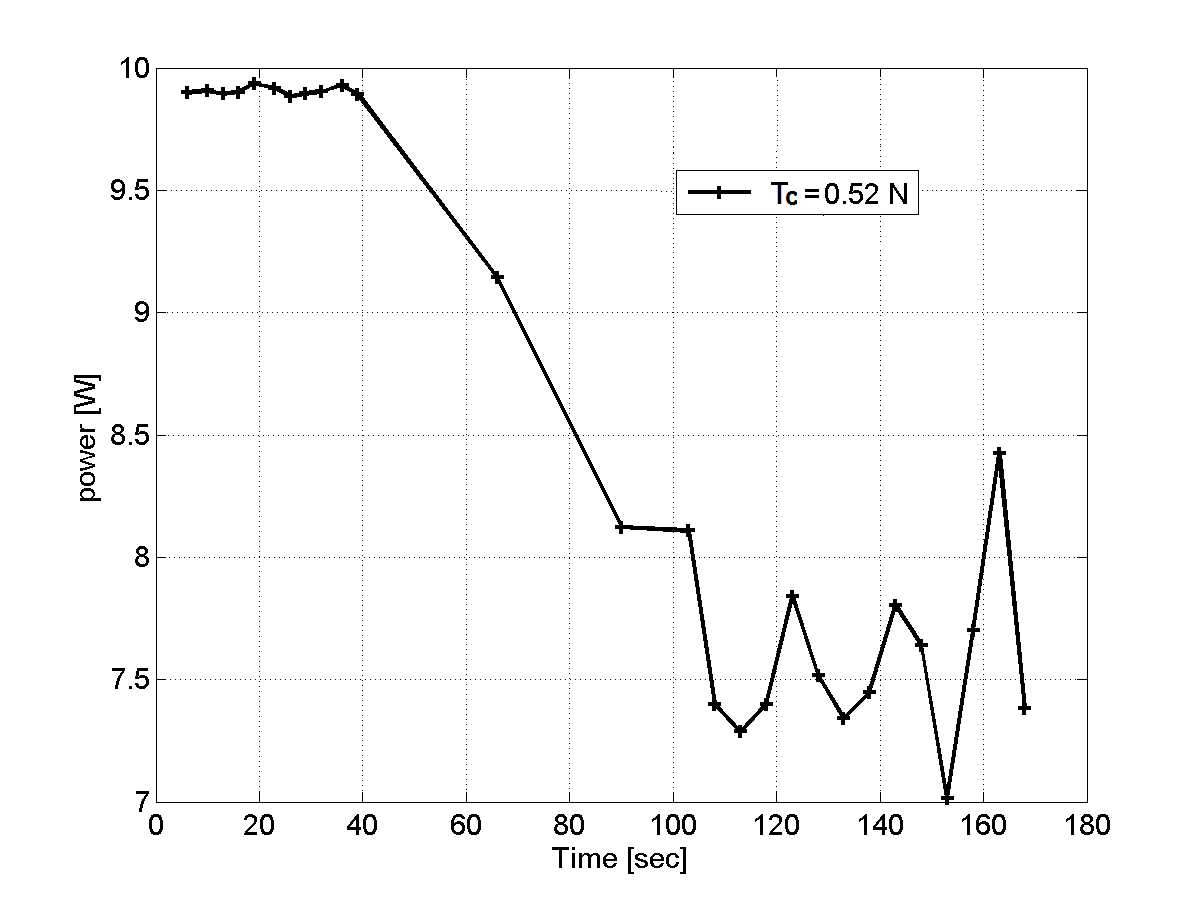}
\caption{Power of the DC Motor during the online optimization process} \label{fig:power_optimization_3}
\end{minipage}
\end{figure}

\begin{figure}[!ht]
\centering
\begin{minipage}[b]{0.48\textwidth}
\centering
\includegraphics[width=0.99\textwidth]{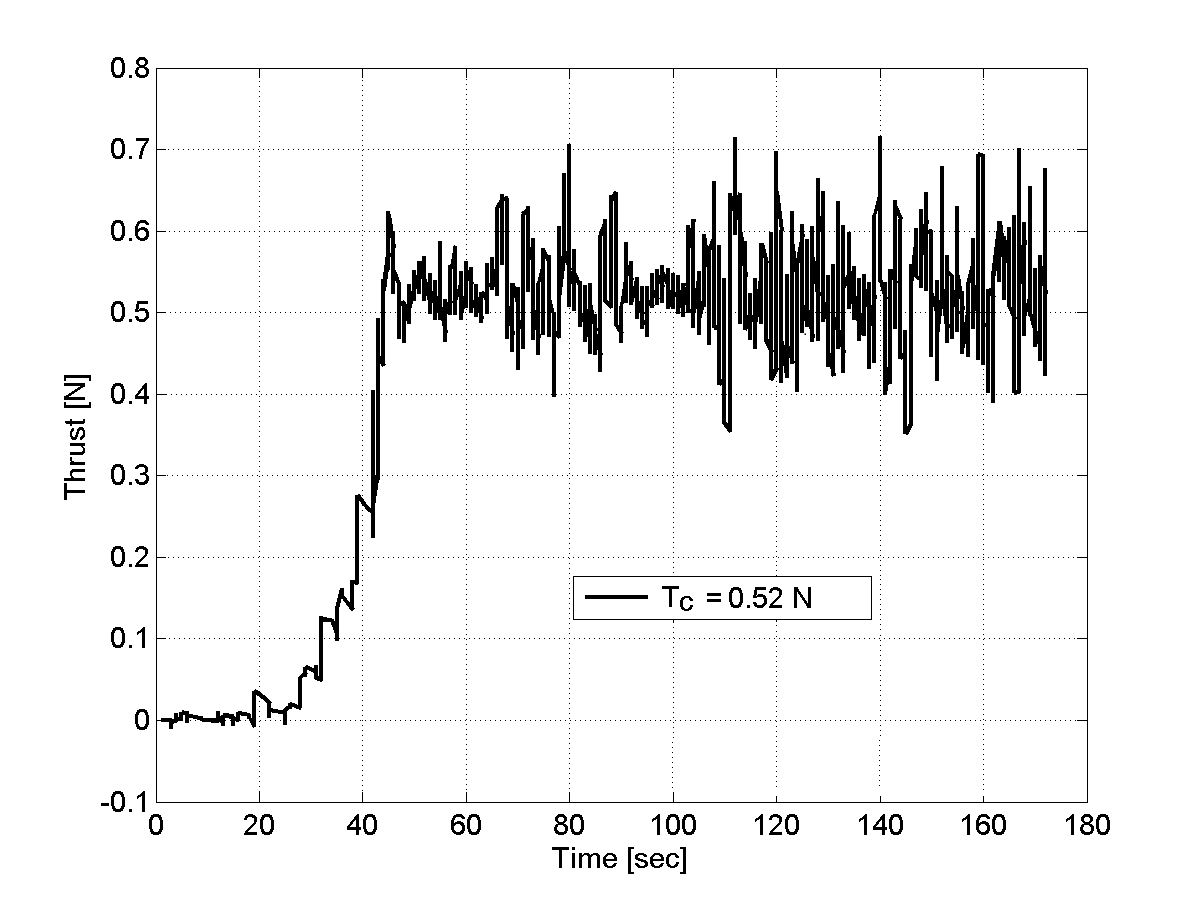}
\caption{Thrust of the propeller during the online optimization process} \label{fig:thrust_optimization_3}
\end{minipage}
\hfill
\begin{minipage}[b]{0.48\textwidth}
\centering
\includegraphics[width=0.99\textwidth]{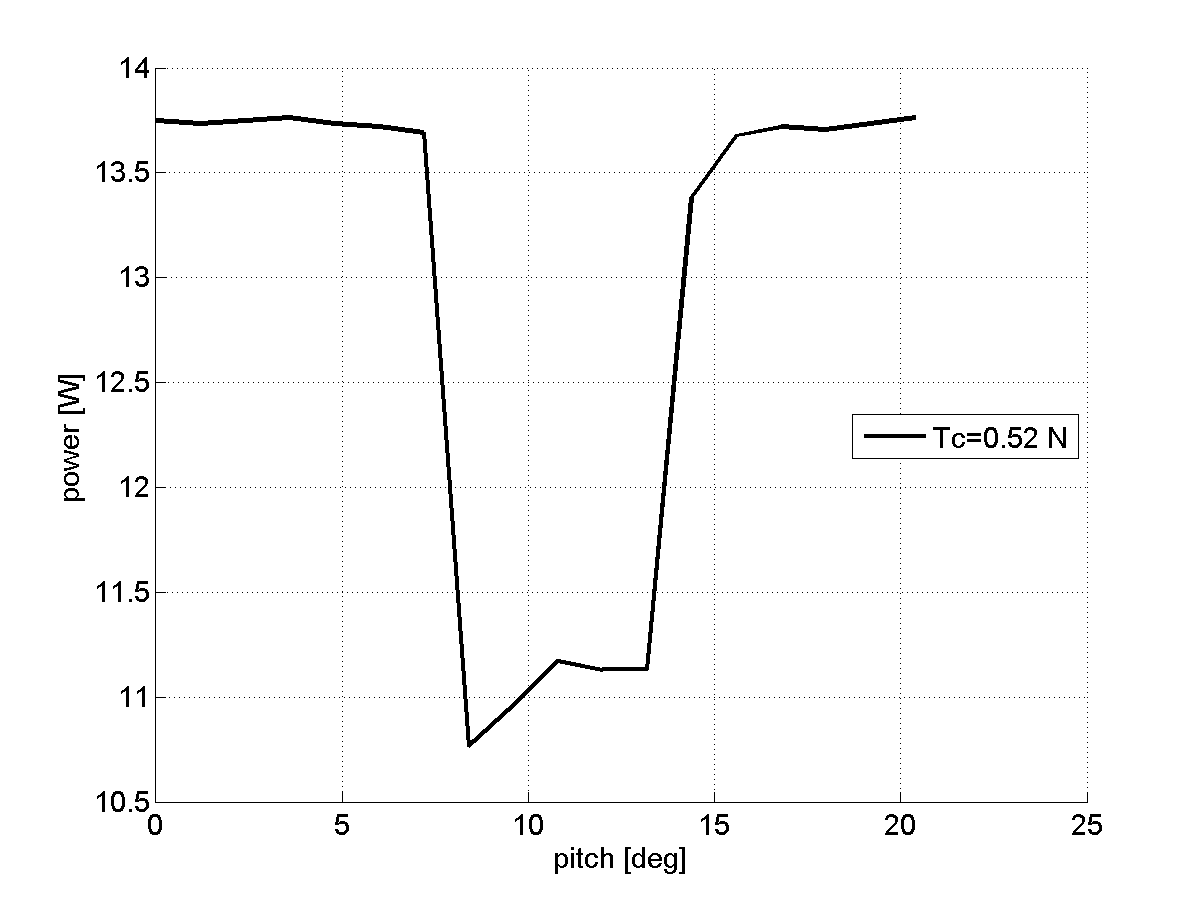}
\caption{Required power versus pitch angle for $T_c=0.52N$} \label{fig:isothrust_minimum}
\end{minipage}
\end{figure}

We now present the experimental results obtained with the basic optimization described above. For the first experiment, the following parameters have been chosen: the thrust command, $T_{c}=0.52~(N)$; the initial pitch angle $\beta_{ini}=0.59~(deg)$; and the incremental step angle $\delta \beta=0.59~(deg)$. Figures \ref{fig:pitch_optimization_3}, \ref{fig:power_optimization_3}, and \ref{fig:thrust_optimization_3} show the evolution of pitch angle, motor power, and propeller thrust as the optimization proceeds.

During the optimization, the pitch angle increases from 0.5 deg to 8 deg approximately, during the first 100 seconds. Then, the pitch angle oscillates between 8.5 deg and 9.5 deg. It can be observed that power decreases from 10 W to 7.4 W during the first 100 seconds, and after that it oscillates between 7.4 W and 7.8 W. It is important to know about the thrust evolution as well, to make sure that the generated thrust of the system is close to the thrust command. As it is apparent in Figure \ref{fig:thrust_optimization_3}, propeller thrust is pretty close to the command; however the system was able to generate the command after 40 seconds only. This is because the thrust command cannot be followed at low pitch angles, and it saturates the intensity available from the power generator before commanded thrust is reached.

The experimental results show that the system finds a stationary point at $\beta = 9$ deg. Figure \ref{fig:isothrust_minimum} shows the evolution of required power versus $\beta$. According to this figure, the basic optimization algorithm has indeed found the minimum power required to track the corresponding thrust (0.52~N). With this optimization approach, it takes too much time for the system to reach minimum power. We now introduce an approach based on variable step length.

\subsection{Improved Optimization: Variable Pitch Steps}%%%%%%%%%%%%%%%%%%%%%%%%%%%%%%%%
The variable step length optimization process is presented as a pseudocode in Algorithm \ref{Alg:VariabPitch}.
\begin{algorithm}[!ht]
 \SetAlgoLined
 {\small
 initialization\;
 $\beta=0.59~deg$;  \newline
 pitch\_step=0.59*3; \newline
 $T_c=0.52~N$; \newline
 prev\_direction=1; \newline
 set\_propeller ($\beta,T_c$); \newline
 measured prev\_power; \newline
 diff\_power=1; \newline
 \While{3 minutes}{
  \eIf{power saturation}{
   next\_direction=1;}{\eIf{diff\_power$>0$}{next\_direction=1;}{next\_direction=-1;}}
   \eIf{prev\_direction $\times$ next\_direction=-1 and step $\ne$ 0.59}{pitch\_step=pitch\_step-0.59}{}
   $\beta=\beta$ + next\_direction $\times$ pitch\_step; \newline
   set\_propeller ($\beta,T_c$); \newline
   measure next\_power; \newline
   diff\_power=prev\_power-next\_power; \newline
   prev\_power=next\_power;     }  }
 \caption{Pseudocode of the improved optimization approach with variable pitch steps}\label{Alg:VariabPitch}
\end{algorithm}

\begin{figure}[!ht]
\centering
\begin{minipage}[b]{0.48\textwidth}
\centering
\includegraphics[width=0.99\textwidth]{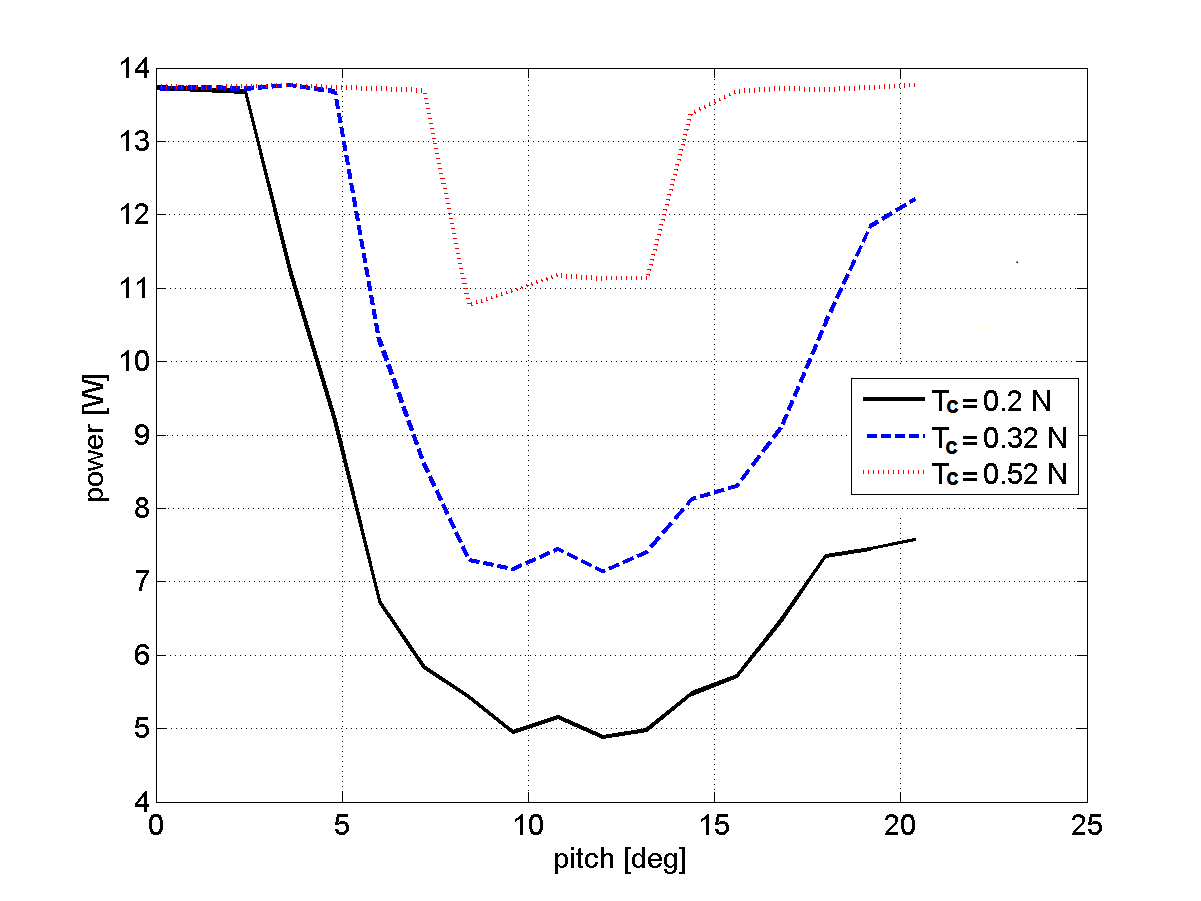}
\caption{Required power versus pitch angle for various values of commanded thrust} \label{fig:Motor_overload}
\end{minipage}
\hfill
\begin{minipage}[b]{0.48\textwidth}
\centering
\includegraphics[width=0.99\textwidth]{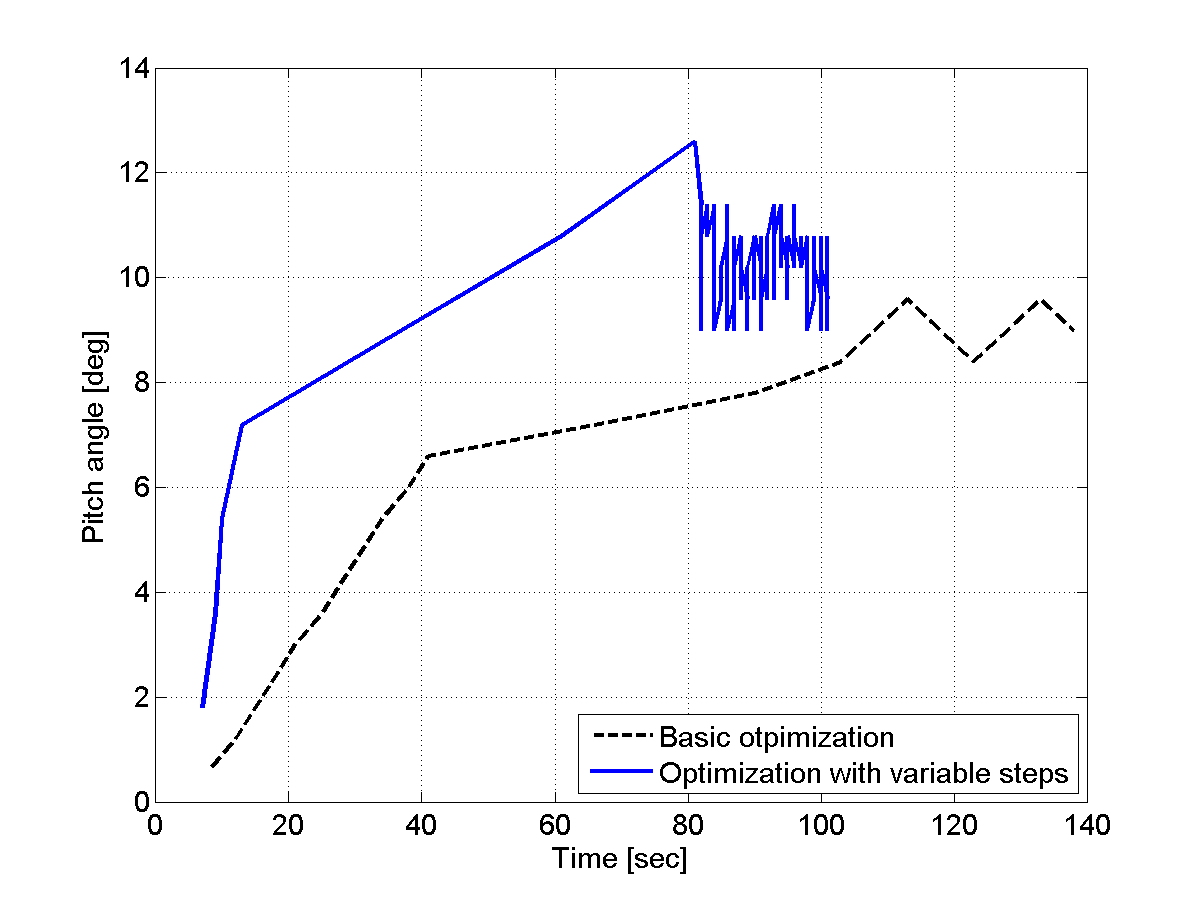}
\caption{Pitch angle evolution for both optimization approaches} \label{fig:pitch_comparaison}
\end{minipage}
\end{figure}

Figure \ref{fig:Motor_overload} shows required power versus pitch angle for various values of commanded thrust. When power saturates near 14 Watt, the system cannot achieve the commanded thrust.

\begin{figure}[!ht]
\centering
\begin{minipage}[b]{0.48\textwidth}
\centering
\includegraphics[width=0.99\textwidth]{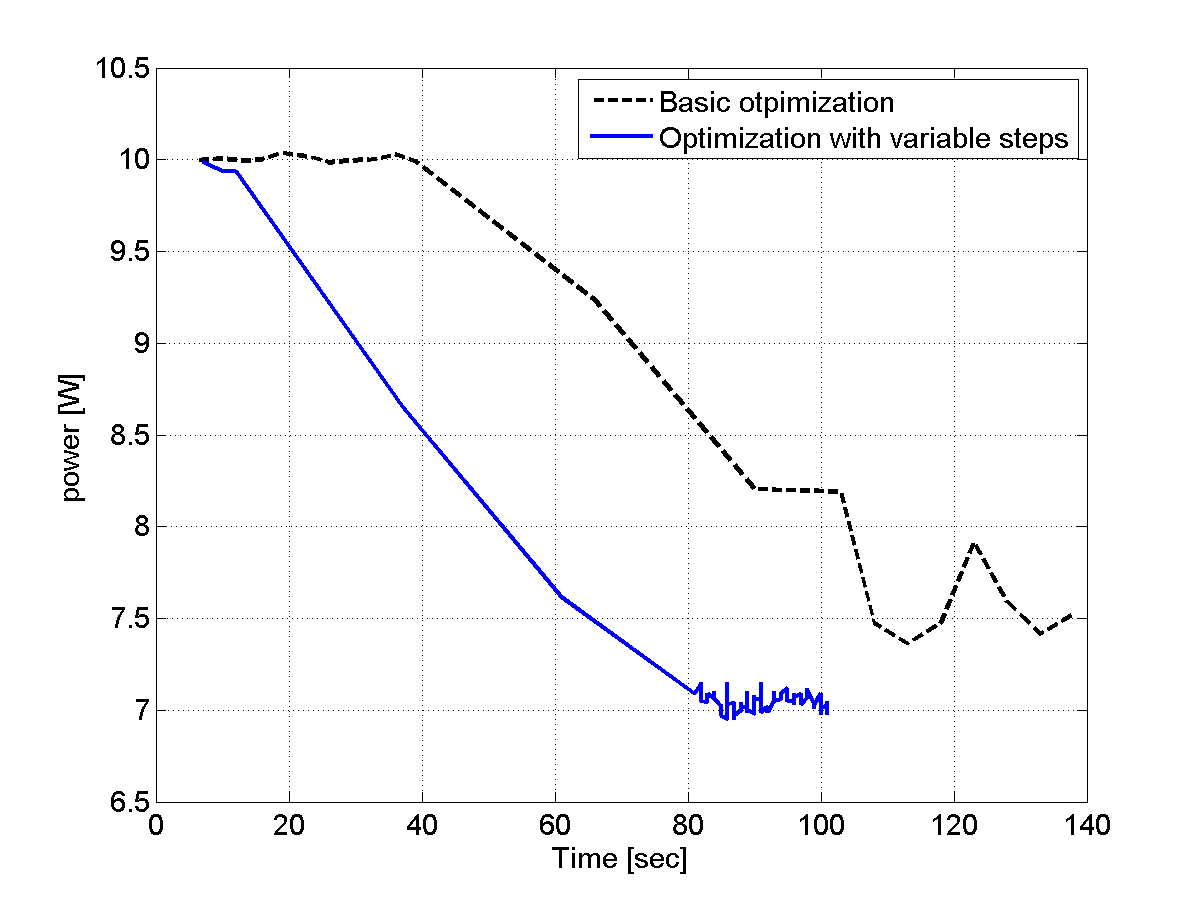}
\caption{Power evolution for both optimization approaches} \label{fig:power_comparaison}
\end{minipage}
\hfill
\begin{minipage}[b]{0.48\textwidth}
\centering
\includegraphics[width=0.99\textwidth]{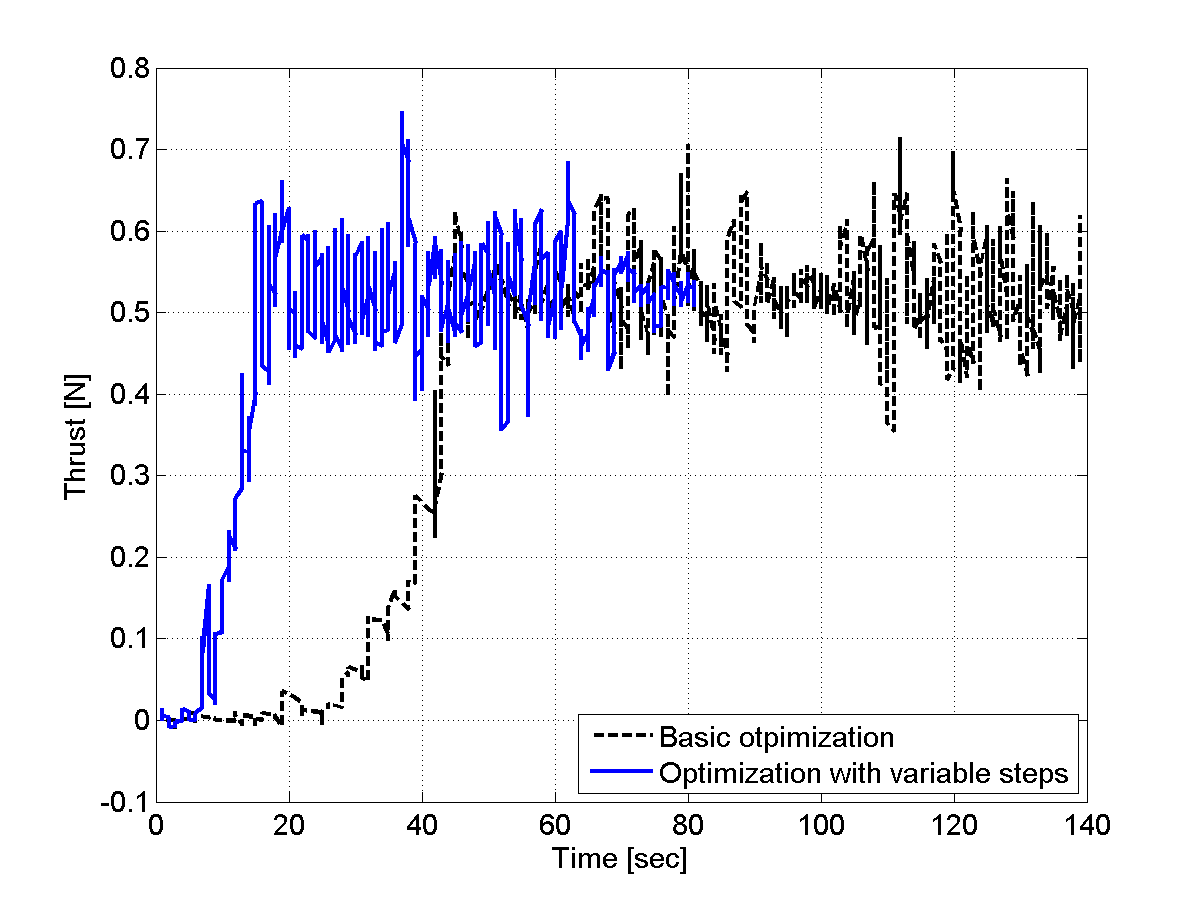}
\caption{Thrust evolution for both optimization approaches} \label{fig:thrust_comparaison}
\end{minipage}
\end{figure}

Figures \ref{fig:pitch_comparaison}, \ref{fig:power_comparaison}, and \ref{fig:thrust_comparaison} show the evolution of the pitch angle, the power, and the thrust for both fixed and variable pitch steps optimization approaches. According to Figures \ref{fig:pitch_comparaison} and \ref{fig:power_comparaison}, the pitch angle converges faster towards the optimal value; the variable pitch steps optimization program requires about 80 seconds, compared with 120 seconds for the fixed pitch steps optimization. In addition, it can be observed that power also decreases much faster, and converges to a slightly better optimum (approximately at about 7 (W)). Figure \ref{fig:thrust_comparaison} shows that considerably less time, about 15 seconds, is needed for the system to reach the commanded thrust value with the variable pitch steps optimization approach, compared with the 45 seconds for the fixed pitch steps optimization approach.

\newpage
%%%%%%%%%%%%%%%%%%%%%%%%%%%%%%%%%%%%%%%%%%%%%%%%%%%%%%%%%%%%%%%%%%%%%%%%%%%%%%%%%%%%%%%%%%%%%%%%%%%%%%%%%%%%%%%%%%%%%%%%%
\section{Conclusion and Future Work}
%%%%%%%%%%%%%%%%%%%%%%%%%%%%%%%%%%%%%%%%%%%%%%%%%%%%%%%%%%%%%%%%%%%%%%%%%%%%%%%%%%%%%%%%%%%%%%%%%%%%%%%%%%%%%%%%%%%%%%%%%
This paper is interested in the online optimization of power for a given variable-pitch propeller propulsion system. It was shown that an online optimization scheme can be implemented to minimize input power for any given, achievable thrust. However, we have shown that online optimization is possible only for sea level atmospheric condition and zero air speed. During flight, conditions are very different than at sea level; thus next steps for our research include testing our system in a wind tunnel at various airspeeds and pressure levels to reflect changes in altitude. Further work includes establishing rigorous proofs of system stability as the online optimization process interacts with real-time thrust control. Such stability proofs are essential in the context of this paper, where safety-critical code must be certified prior to routine operations.

%%%%%%%%%%%%%%%%%%%%%%%%%%%%%%%%%%%%%%%%%%%%%%%%%%%%%%%%%%%%%%%%%%%%%%%%%%%%%%%%
\section*{Acknowledgment}
This material is based upon the work supported by the National Aeronautics and Space Administration (NASA), the National Science Foundation (NSF), the Army Research Office, and also the Agence Nationale de la Recherche/DGA - ASTRID program.
%%%%%%%%%%%%%%%%%%%%%%%%%%%%%%%%%%%%%%%%%%%%%%%%%%%%%%%%%%%%%%%%%%%%%%%%%%%%%%%%

\bibliographystyle{plain}
%\bibliography{DistGainSchedulTurbineControl}

%%%%%%%%%%%%%%%%%%%%%%%%%%%%%%%%%%%%%%%%%%%%%%%%%%%%%%%%%%%%%%%%%%%%%%%
%\appendix       %%% starting appendix
%\section*{Appendix A: Head of First Appendix}
%Avoid Appendices if possible.

%%%%%%%%%%%%%%%%%%%%%%%%%%%%%%%%%%%%%%%%%%%%%%%%%%%%%%%%%%%%%%%%%%%%%%
%\end{multicols}

\end{document}